\input amsppt.sty
\magnification=\magstep1
\hsize=30truecc
\vsize=22.2truecm
\baselineskip=16truept
\nologo
\pageno=1
\TagsOnRight
\topmatter

\def\C{\Bbb C}

\def\l{\left}
\def\r{\right}
\def\bg{\bigg}
\def\({\bg(}
\def\[{\bg[}
\def\){\bg)}
\def\]{\bg]}
\def\t{\text}
\def\f{\frac}

\def\bi{\binom}

\def\cs{\ldots}
\def\ls{\leqslant}
\def\gs{\geqslant}

\def\Da{\Delta}

\def\bi{\binom}

\def\Proof{\noindent{\it Proof}}

\topmatter
\hbox{Preprint (August 5, 2004), {\tt arXiv:math.NT/0407363}.}
\bigskip
\title {New identities involving Bernoulli and Euler polynomials}\endtitle
\rightheadtext{Identities involving Bernoulli and Euler polynomials}
\author {Hao Pan and Zhi-Wei Sun}\endauthor
\abstract In this paper we obtain several new identities for Bernoulli
and Euler polynomials; some of them extend Miki's and Matiyasevich's identities.
Our new method involves differences and derivatives of polynomials.
\endabstract
\thanks  2000 {\it Mathematics Subject Classification}.
Primary 11B68; Secondary 05A19.
\newline\indent The second author is responsible for all the communications
and supported by
the National Natural Science Foundation of P. R. China.
\endthanks
\endtopmatter
\document
\hsize=30truecc
\vsize=22.2truecm
\baselineskip=16truept

\heading 1. Introduction\endheading The Bernoulli numbers $B_0,
B_1,\ldots$ are given by
$$\frac{z}{e^z-1}=\sum_{n=0}^{\infty}B_n\frac{z^n}{n!}\qquad (|z|<2\pi),$$
they play important roles in many aspects.
Most research on Bernoulli numbers concentrates on their congruence properties (cf. e.g. [Su1]).
However, there are also some interesting identities concerning Bernoulli numbers (see, e.g. [Di] and [Su2]).

In 1978 Miki [Mi] proposed the following curious identity
which involves both an ordinary convolution and a binomial convolution
of Bernoulli numbers:
$$\sum_{k=2}^{n-2}\frac{B_kB_{n-k}}{k(n-k)}-\sum_{l=2}^{n-2}\bi{n}{l}\frac{B_lB_{n-l}}{l(n-l)}
=\frac{2}nH_nB_n\tag 1.1$$
for any $n=4,5,\ldots$, where
$$H_n=1+\frac{1}{2}+\cdots+\frac{1}{n}.$$
In the original proof of this identity,
Miki showed that the two sides of $(1.1)$ are congruent modulo all sufficiently large primes.
Shiratani and Yokoyama [SY] gave another proof of $(1.1)$ by $p$-adic analysis, and
 Gessel [Ge] reproved Miki's identity (1.1)
 by using the ordinary generating function and the exponential
generating function of Stirling numbers of the second kind.

Inspired by Miki's work, Matiyasevich [Ma] found the following two
identities of the same nature by the software {\it Mathematica}.
$$\sum_{k=2}^{n-2}\frac{B_k}kB_{n-k}-\sum_{l=2}^{n-2}\bi{n}{l}\frac{B_l}lB_{n-l}=H_nB_n\tag 1.2$$
and $$(n+2)\sum_{k=2}^{n-2}B_kB_{n-k}-2\sum_{l=2}^{n-2}\bi{n+2}{l}B_lB_{n-l}=n(n+1)B_n
\tag 1.3$$
for each $n=4,5,\ldots$.
We mention that (1.2) is actually equivalent to Miki's identity (1.1).
The reason is as follows:
$$\align&\sum_{k=2}^{n-2}\f{B_kB_{n-k}}{k(n-k)}-\sum_{l=2}^{n-2}\bi nl\f{B_lB_{n-l}}{l(n-l)}
\\=&\f1n\sum_{k=2}^{n-2}\(\f1k+\f1{n-k}\)B_kB_{n-k}
-\f1n\sum_{l=2}^{n-2}\bi{n}{l}\(\f1l+\f1{n-l}\)B_lB_{n-l}
\\=&\f2n\sum_{k=2}^{n-2}\f{B_k}kB_{n-k}
-\f2n\sum_{l=2}^{n-2}\bi nl\f{B_l}lB_{n-l}.
\endalign$$

Quite recently Dunne and Schubert [DS] presented a new approach to (1.1) and (1.3)
motivated by quantum field theory and string theory.

The Bernoulli polynomials $B_n(x)=\sum_{k=0}^n\bi{n}{k}B_{k}x^{n-k}$ $(n=0,1,2,\ldots)$
have the following basic properties:
$$B_n(x+1)-B_n(x)=nx^{n-1},\  B_n(x+y)=\sum_{k=0}^n\bi{n}{k}B_{k}(x)y^{n-k},$$
and also $B_n'(x)=nB_{n-1}(x)$ for $n>0$.

In this paper we extend Miki's identity (1.1) and Matiyasevich's
identity (1.3) to Bernoulli polynomials.

\proclaim{Theorem 1.1} Let $n>1$ be an integer. Then
$$\aligned
&\sum_{k=1}^{n-1}\frac{B_k(x)B_{n-k}(y)}{k(n-k)}
-\sum_{l=1}^{n}\bi{n-1}{l-1}\frac{B_l(x-y)B_{n-l}(y)+B_l(y-x)B_{n-l}(x)}{l^2}\\
&\qquad=\frac{H_{n-1}}{n}(B_n(x)+B_n(y))+\frac{B_n(x)-B_n(y)}{n(x-y)}
\endaligned\tag 1.4$$
and
$$\aligned
&\sum_{k=0}^nB_k(x)B_{n-k}(y)-\sum_{l=0}^n\bi{n+1}{l+1}\frac{B_l(x-y)B_{n-l}(y)+B_l(y-x)B_{n-l}(x)}{l+2}\\
&\qquad=\frac{B_{n+1}(x)+B_{n+1}(y)}{(x-y)^2}-\f{2}{n+2}\cdot\frac{B_{n+2}(x)-B_{n+2}(y)}{(x-y)^3}.
\endaligned\tag 1.5$$
\endproclaim
We remark that (1.4) has the following equivalent version
$$\aligned
&2\sum_{k=1}^{n-1}\frac{B_k(x)}{k}B_{n-k}(y)
-H_{n-1}(B_n(x)+B_n(y))-\frac{B_n(x)-B_n(y)}{x-y}\\
&\qquad=\sum_{l=1}^{n}\bi{n}{l}\(\frac{B_l(x-y)}{l}B_{n-l}(y)+\frac{B_l(y-x)}{l}B_{n-l}(x)\).
\endaligned\tag 1.4$'$ $$

\proclaim{Corollary 1.1} Let $n\gs2$ be an integer. Then we have
$$\sum_{k=1}^{n-1}\frac{B_k(x)B_{n-k}(x)}{k(n-k)}-2\sum_{l=2}^n\bi{n-1}{l-1}\frac{B_lB_{n-l}(x)}{l^2}
=\frac{2}{n}H_{n-1}B_n(x),
\tag 1.6$$
and
$$\sum_{k=0}^nB_k(x)B_{n-k}(x)-2\sum_{l=2}^n\bi{n+1}{l+1}\frac{B_lB_{n-l}(x)}{l+2}=(n+1)B_n(x).\tag1.7$$
\endproclaim
\Proof. Letting $y$ tend to $x$ and recalling that
$B_n'(x)=nB_{n-1}(x)$, we immediately find that (1.4) implies
(1.6).

Now we go to prove (1.7).  Let $P(z)=B_{n+2}(z)/(n+2)$. Then
$P'(z)=B_{n+1}(z)$, $P''(z)=(n+1)B_n(z)$ and
$P'''(z)=n(n+1)B_{n-1}(z)$. In light of Taylor's expansion,
$$P(y)-P(x)=P'(x)(y-x)+\frac{P''(x)}{2!}(y-x)^2+\frac{P'''(x)}{3!}(y-x)^3+\cdots$$
and
$$P'(y)-P'(x)=P''(x)(y-x)+\frac{P'''(x)}{2!}(y-x)^2+\cdots.$$
Therefore
$$\align
&\lim_{y\to x}\(\frac{B_{n+1}(x)+B_{n+1}(y)}{(x-y)^2}
-\f 2{n+2}\cdot\frac{B_{n+2}(x)-B_{n+2}(y)}{(x-y)^3}\)\\
=&\lim_{y\to x}\(\f{P'(x)+P'(y)}{(x-y)^2}-\frac{2(P(x)-P(y))}{(x-y)^3}\)\\
=&\lim_{y\to x}\(\f{P'(y)-P'(x)}{(y-x)^2}-2\frac{P(y)-P(x)-P'(x)(y-x)}{(y-x)^3}\)\\
=&\lim_{y\to x}\(\f{P''(x)}{y-x}+\frac{P'''(x)}{2!}+\cdots-2\(\frac{P''(x)}{2!(y-x)}
+\frac{P'''(x)}{3!}+\cdots\)\)\\
=&\f{P'''(x)}6=\frac{n(n+1)}{6}B_{n-1}(x).
\endalign$$
In view of this, we can easily get (1.7) from (1.5) by letting $y$
tend to $x$. \qed

Now let us see how Miki's identity follows from (1.6). In fact, (1.6) in the case $x=0$ yields that
$$\align
\sum_{k=1}^{n-1}\frac{B_kB_{n-k}}{k(n-k)}=&2\sum_{l=1}^{n-1}\bi{n-1}{l}\frac{B_lB_{n-l}}{l(n-l)}
+\frac{2B_n}{n^2}-2B_1B_{n-1}+\frac{2}nH_{n-1}B_n\\
=&\sum_{l=1}^{n-1}\(\bi{n-1}{l}+\bi{n-1}{n-l}\)\frac{B_lB_{n-l}}{l(n-l)}
+\frac2n H_nB_n+B_{n-1}\\
=&\sum_{l=1}^{n-1}\bi{n}{l}\frac{B_lB_{n-l}}{l(n-l)}+\frac{2}{n}H_nB_n+B_{n-1}.
\endalign$$
We also mention that (1.7) in the case $x=0$ gives Matiyasevich's
identity (1.3).

\proclaim{Corollary 1.2} Let $n\gs4$ be an integer. Then
$$\sum_{k=2}^{n-2}\f{\bar B_k}k\bar B_{n-k}=\f n2\sum_{k=2}^{n-2}\frac{\bar B_{k}\bar B_{n-k}}{k(n-k)}
=\sum_{k=2}^{n}\bi{n}{k}\frac{B_{k}}{k}\bar B_{n-k}+H_{n-1}\bar B_n,$$
where $\bar B_{k}=(2^{1-k}-1)B_k$.
\endproclaim
\Proof. Simply take  $x=1/2$ in $(1.6)$ and use the known formula
$B_n(1/2)=\bar B_n$. (Note also that
$n/(k(n-k))=1/k+1/(n-k)$.)\qed

 The last equality in Corollary 1.2 was first found by C. Faber and R. Pandharipande,
and then confirmed by Zagier (cf. [FP]).

The Euler polynomials $E_n(x)\ (n=0,1,2,\ldots)$ are defined by
$$\frac{2e^{xz}}{e^z+1}=\sum_{n=0}^{\infty}E_n(x)\frac{z^n}{n!}.$$
Here are some basic properties of Euler polynomials:
$$E_n(x+1)+E_n(x)=2x^n,\ E_n(x+y)=\sum_{k=0}^n\bi{n}{k}E_k(x)y^{n-k},$$
and also $E_n'(x)=nE_{n-1}(x)$ if $n>0$. It is also known that
$$E_n(x)=\f 2{n+1}\(B_{n+1}(x)-2^{n+1}B_{n+1}\l(\f x2\r)\).$$

  Similar to Theorem 1.1 we have the following identities involving Euler polynomials.
\proclaim{Theorem 1.2} Let $n$ be a positive integer. Then
$$\aligned
&\sum_{k=0}^nE_k(x)E_{n-k}(y)-\frac{4}{n+2}\cdot\f{B_{n+2}(x)-B_{n+2}(y)}{x-y}
\\=&-2\sum_{l=0}^{n+1}\bi{n+1}{l}\frac{E_l(x-y)B_{n+1-l}(y)+E_l(y-x)B_{n+1-l}(x)}{l+1}.
\endaligned\tag 1.8$$
Also,
$$\aligned &\sum_{k=1}^n\frac{B_k(x)}kE_{n-k}(y)-H_nE_n(y)-\frac{E_n(x)-E_n(y)}{x-y}
\\=&\sum_{l=1}^n\bi{n}{l}\(\f{B_l(x-y)}lE_{n-l}(y)-\f{E_{l-1}(y-x)}2E_{n-l}(x)\),
\endaligned\tag 1.9$$
and
$$\aligned
&\sum_{k=0}^nB_k(x)E_{n-k}(y)
\\=&\sum_{l=1}^n\bi{n+1}{l+1}\(B_l(x-y)E_{n-l}(y)-\frac{E_{l-1}(y-x)}2E_{n-l}(x)\)
\\&+(n+1)\(\frac{E_n(x)}{x-y}+E_n(y)\)-\frac{E_{n+1}(x)-E_{n+1}(y)}{(x-y)^2}.
\endaligned\tag 1.10$$
\endproclaim

\proclaim{Corollary 1.3} Let $n$ be any nonnegative integer. Then we have
$$\gather
(n+2)\sum_{k=0}^nE_k(x)E_{n-k}(x)=8\sum_{l=2}^{n+2}\bi{n+2}{l}(2^{l}-1)\f{B_l}lB_{n+2-l}(x),\tag1.11
\\\sum_{k=1}^n\frac{B_k(x)}kE_{n-k}(x)-\sum_{l=2}^n\bi{n}{l}2^l\f{B_l}lE_{n-l}(x)=H_nE_n(x),\tag1.12
\\\sum_{k=0}^nB_k(x)E_{n-k}(x)-\sum_{l=2}^n\bi{n+1}{l+1}(2^l+l-1)\f{B_l}lE_{n-l}(x)=(n+1)E_n(x).\tag1.13
\endgather$$
\endproclaim
\Proof. Letting $y$ tend to $x$
and noting that $E_l(0)=2(1-2^{l+1})B_{l+1}/(l+1)$,
we then obtain the (1.11) and (1.12) from (1.8) and (1.9) respectively.

 Since
 $$E_{n+1}(y)-E_{n+1}(x)=E_{n+1}'(x)(y-x)+\f{E_{n+1}''(x)}{2!}(y-x)^2+\cdots,$$
we have
$$\align
&\lim_{y\to x}\((n+1)\f{E_n(x)}{x-y}-\f{E_{n+1}(x)-E_{n+1}(y)}{(x-y)^2}\)\\
=&\lim_{y\to x}\f{E_{n+1}(y)-E_{n+1}(x)-(y-x)E_{n+1}'(x)}{(y-x)^2}\\
=&\f{E_{n+1}''(x)}{2!}=\frac{n(n+1)}{2}E_{n-1}(x).
\endalign$$
Thus, (1.13) follows from (1.10) by letting $y$ tend to $x$. We are done. \qed

In the next section we will prove Theorems 1.1 and 1.2 by a new approach via
differences and derivatives of polynomials.
In Section 3 we will give a further extension of Corollary 1.1 with help of the gamma function.

\heading 2. Proofs of Theorems 1.1 and 1.2\endheading

\proclaim{Lemma 2.1} Let $P(x),Q(x)\in \C[x]$ where $\C$ is the field of complex numbers.

{\rm (i)} We have
$$\aligned\Delta(P(x)Q(x))=&P(x)\Delta(Q(x))+Q(x)\Delta(P(x))+\Delta(P(x))\Delta(Q(x)),
\\=&\Delta^*(P(x))\Delta^*(Q(x))-P(x)\Delta^*(Q(x))-Q(x)\Delta^*(P(x))
\endaligned$$
and
$$\Delta^*(P(x)Q(x))=\Delta(P(x))\Delta^*(Q(x))+P(x)\Delta^*(Q(x))-Q(x)\Delta(P(x)),$$
where the operators $\Delta$ and $\Delta^*$ are given by $\Delta(f(x))=f(x+1)-f(x)$
and $\Delta^*(f(x))=f(x+1)+f(x)$.

{\rm (ii)} If $\Delta(P(x))=\Delta(Q(x))$, then $P'(x)=Q'(x)$.
 If $\Delta^*(P(x))=\Delta^*(Q(x))$, then $P(x)=Q(x)$.
\endproclaim
\Proof. (i) Part (i) can be verified directly.

(ii) Suppose that $\Delta(P(x))=\Delta(Q(x))$. Then,
$$P(n)-P(0)=\sum_{k=0}^{n-1}\Delta(P(k))=\sum_{i=0}^{n-1}\Delta(Q(k))=Q(n)-Q(0)$$
for every $n=1,2,3,\ldots$.
Now that the polynomial $g(x)=P(x)-Q(x)-P(0)+Q(0)$ has infinitely many zeroes,
we must have $g(x)=0$ and hence $P'(x)=Q'(x)$.

Now assume that $\Delta^*(P(x))=\Delta^*(Q(x))$. Then
$$P(n)-Q(n)=-(P(n-1)-Q(n-1))=\cs=(-1)^n(P(0)-Q(0))$$
for every $n=1,2,\ldots$.
Since the equations $P(x)-Q(x)=P(0)-Q(0)$ and $P(x)-Q(x)=-(P(0)-Q(0))$ both have infinitely many roots,
we must have $P(x)=Q(x)$.

 The proof of Lemma 2.1 is now complete.
\qed

\proclaim{Lemma 2.2} Let
$n$ be any positive integer. Then
$$\sum_{k=1}^n\f{B_k(x+y)}kx^{n-k}=\sum_{l=1}^n\f{B_l(y)}lx^{n-l}+H_nx^n\tag2.1$$
and
$$\sum_{k=0}^nE_k(x+y)x^{n-k}=\sum_{l=0}^n\bi{n+1}{l+1}E_l(y)x^{n-l}.\tag2.2$$
\endproclaim
\Proof. Note that
$$\align\sum_{k=1}^n\f{B_k(x+y)}kx^{n-k}
=&\sum_{k=1}^n\f{x^{n-k}}k\(\sum_{l=1}^k\bi klB_l(y)x^{k-l}+x^k\)
\\=&\sum_{l=1}^n\f{B_l(y)}lx^{n-l}\sum_{k=l}^n\bi{k-1}{l-1}+\sum_{k=1}^n\f{x^n}k
\\=&\sum_{l=1}^n\bi nl\f{B_l(y)}lx^{n-l}+H_nx^n
\endalign$$
where in the last step we apply a well-known identity of Chu (see, e.g. [GKP, (5.10)]).
Similarly, we have
$$\align&\sum_{k=0}^{n}E_k(x+y)x^{n-k}=\sum_{k=0}^{n}x^{n-k}\sum_{l=0}^k\bi{k}{l}E_l(y)x^{k-l}
\\=&\sum_{l=0}^{n}E_l(y)x^{n-l}\sum_{k=l}^{n}\bi{k}{l}=\sum_{l=0}^{n}\bi{n+1}{l+1}E_l(y)x^{n-l}.
\endalign$$
So both (2.1) and (2.2) hold. \qed

\medskip
{\noindent{\it Proof of Theorem 1.1}}. Observe that
$$\align
&\frac{\partial}{\partial x}\frac{\partial}{\partial y}\sum_{l=1}^{n}\bi{n-1}{l-1}\frac{B_l(y-x)B_{n-l}(x)}{l^2}
\\=&\f{\partial}{\partial x}\sum_{l=1}^n\bi{n-1}{l-1}\f{B_{l-1}(y-x)}lB_{n-l}(x)
\\=&\sum_{l=1}^{n-1}\bi{n-1}{l-1}\f{n-l}lB_{l-1}(y-x)B_{n-1-l}(x)
\\&-\sum_{l=2}^n\bi{n-1}{l-1}\f{l-1}lB_{l-2}(y-x)B_{n-l}(x)
\\=&\sum_{l=1}^{n-1}\bi{n-1}lB_{l-1}(y-x)B_{n-1-l}(x)
\\&+\sum_{l=2}^n\bi{n-1}{l-1}\l(\f1l-1\r)B_{l-2}(y-x)B_{n-l}(x)
\\=&\sum_{l=0}^{n-2}\bi{n-1}{l+1}\frac{B_{l}(y-x)B_{n-2-l}(x)}{l+2}.
\endalign$$
Therefore
$$\align
&\frac{\partial}{\partial x}\frac{\partial}{\partial y}\sum_{l=1}^{n}\bi{n-1}{l-1}\frac{B_l(x-y)B_{n-l}(y)}{l^2}
\\=&\frac{\partial}{\partial y}\frac{\partial}{\partial x}\sum_{l=1}^{n}\bi{n-1}{l-1}\frac{B_l(x-y)B_{n-l}(y)}{l^2}
\\=&\sum_{l=0}^{n-2}\bi{n-1}{l+1}\frac{B_{l}(x-y)B_{n-2-l}(y)}{l+2}.
\endalign$$
We also have
$$\align
&\f{\partial}{\partial x}\f{\partial}{\partial y}\(\frac{B_n(y)-B_n(x)}{n(y-x)}\)\\
=&\f{\partial}{\partial x}\(\frac{nB_{n-1}(y)}{n(y-x)}-\frac{B_n(y)-B_n(x)}{n(y-x)^2}\)\\
=&\frac{B_{n-1}(y)}{(y-x)^2}+\frac{B_{n-1}(x)}{(y-x)^2}-\f2n\cdot\frac{B_n(x)-B_n(y)}{(x-y)^3}.
\endalign$$
Let $L(x,y)$ and $R(x,y)$ denote the left hand side and right hand side of (1.4) respectively.
In view of the above, that
$\f{\partial}{\partial x}\f{\partial}{\partial y}L(x,y)=\f{\partial}{\partial x}\f{\partial}{\partial y}R(x,y)$
gives (1.5) with $n$ replaced by $n-2$.

If we substitute $x+y$ for $y$ in (1.4), then we get the following
equivalent version of (1.4).
$$\aligned
&\sum_{k=1}^{n-1}\frac{B_k(x+y)B_{n-k}(x)}{k(n-k)}\\
=&\sum_{l=1}^{n}\f1{l^2}\bi{n-1}{l-1}\l(B_l(y)B_{n-l}(x)+B_l(-y)B_{n-l}(x+y)\r)\\
&+\frac{H_{n-1}}n(B_n(x+y)+B_n(x))+\frac{B_n(x+y)-B_n(x)}{ny}.
\endaligned\tag2.3$$

Now it suffices to prove (2.3) only. Let us view $y$ as a fixed parameter.
Denote by $P_n(x)$ and $Q_n(x)$ the left hand side and the right hand side of (2.3) respectively.
It is easy to check that
$$P_{n+1}'(x)-nP_n(x)=\f{B_n(x+y)+B_n(x)}n=Q_{n+1}'(x)-nQ_n(x).$$
By Lemma 2.1(ii) we need only to show $\Delta(P_{n+1}(x))=\Delta(Q_{n+1}(x))$.

In view of Lemma 2.1(i) and the fact that $\Da(B_k(x))=kx^{k-1}$,
$$\aligned &\Delta(P_{n+1}(x))=\sum_{k=1}^n\Delta\(\frac{B_k(x+y)}k\cdot\f{B_{n+1-k}(x)}{n+1-k}\)\\
=&\sum_{k=1}^n\frac{B_k(x+y)}{k}x^{n-k}
+\sum_{k=1}^{n}\frac{B_{n+1-k}(x)}{n+1-k}(x+y)^{k-1}
+\sum_{k=1}^{n}(x+y)^{k-1}x^{n-k}
\\=&\sum_{k=1}^n\f{B_k(x+y)}kx^{n-k}+\sum_{k=1}^n\f{B_k(x)}k(x+y)^{n-k}+\sum_{k=0}^{n-1}(x+y)^kx^{n-1-k}.
\endaligned$$
Applying Lemma 2.2 we then get
$$\align
\Delta(P_{n+1}(x))=&
\sum_{l=1}^n\bi{n}{l}\frac{B_l(y)}lx^{n-l}+H_nx^n
+\sum_{l=1}^n\bi{n}{l}\frac{B_l(-y)}l(x+y)^{n-l}
\\&+H_n(x+y)^n+\frac{(x+y)^n-x^n}{(x+y)-x}.
\endalign$$
On the other hand, it is easy to see that $\Delta(Q_{n+1}(x))$ also equals the right hand side of
the last equality.
Therefore $\Delta(P_{n+1}(x))=\Delta(Q_{n+1}(x))$ as required.
This concludes our proof. \qed

\medskip
\noindent{\it Proof of Theorem 1.2}.
Substituting $x+y$ for $y$ in (1.8) we then get the following equivalent form of (1.8):
$$\aligned
&\sum_{k=0}^nE_k(x+y)E_{n-k}(x)-\frac{4}{n+2}\cdot\f{B_{n+2}(x+y)-B_{n+2}(x)}{y}
\\=&-2\sum_{l=0}^{n+1}\bi{n+1}{l}\frac{E_l(y)B_{n+1-l}(x)+E_l(-y)B_{n+1-l}(x+y)}{l+1}.
\endaligned\tag2.4$$
If we substitute $x+y$ for $x$ and $x$ for $y$ in (1.9),
we then have the following equivalent version of (1.9):
$$\aligned &\sum_{k=1}^n\frac{B_k(x+y)}kE_{n-k}(x)-H_nE_n(x)-\frac{E_n(x+y)-E_n(x)}{y}
\\=&\sum_{l=1}^n\bi{n}{l}\(\frac{B_l(y)}lE_{n-l}(x)-\f{E_{l-1}(-y)}2E_{n-l}(x+y)\).
\endaligned\tag 2.5$$
Note that (1.10) with $n$ replaced by $n-1$
follows from (1.9) by taking partial derivatives with respect to $x$.
In view of the above, we only need to prove (2.4) and (2.5) with $y$ fixed.

(a) For $0\ls k\ls n$, by Lemma 2.1(i) we have
$$\align&\Delta(E_k(x+y)E_{n-k}(x))
\\=&\Delta^*(E_k(x+y))\Delta^*(E_{n-k}(x))
\\&-E_k(x+y)\Delta^*(E_{n-k}(x))-E_{n-k}(x)\Delta^*(E_k(x+y))
\\=&2(x+y)^k\cdot 2x^{n-k}-E_k(x+y)\cdot2x^{n-k}-E_{n-k}(x)\cdot 2(x+y)^k.
\endalign$$
Thus
$$\align&\Delta\(\sum_{k=0}^{n}E_k(x+y)E_{n-k}(x)\)
-4\sum_{k=0}^{n}(x+y)^kx^{n-k}
\\=&-2\sum_{k=0}^{n}E_k(x+y)x^{n-k}-2\sum_{k=0}^{n}E_{n-k}(x)(x+y)^k.
\endalign$$
With help of Lemma 2.2, we have
$$\align&\Delta\(\sum_{k=0}^{n}E_k(x+y)E_{n-k}(x)\)-4\f{(x+y)^{n+1}-x^{n+1}}{(x+y)-x}
\\&=-2\sum_{l=0}^{n}\bi{n+1}{l+1}E_l(y)x^{n-l}-2\sum_{l=0}^{n}\bi{n+1}{l+1}E_l(-y)(x+y)^{n-l}.
\endalign$$
From this we can easily check that $\Delta(P_n(x))=\Delta(Q_n(x))$
where $P_n(x)$ and $Q_n(x)$ denote the left hand side and the right hand side of
(2.4) respectively.

Now that $\Delta(P_{n+1}(x))=\Delta(Q_{n+1}(x))$, we have $P'_{n+1}(x)=Q'_{n+1}(x)$ by Lemma 2.1(ii).
Clearly $P'_{n+1}(x)=(n+2)P_n(x)$ and $Q'_{n+1}(x)=(n+2)Q_n(x)$. So $P_n(x)=Q_n(x)$ and hence (2.4) holds.

(b) Now we turn to prove (2.5). For $1\ls k\ls n$, by Lemma 2.1(i) we have
$$\align&\Delta^*(B_k(x+y)E_{n-k}(x))
\\=&\Delta(B_k(x+y))\Delta^*(E_{n-k}(x))
\\&+B_k(x+y)\Delta^*(E_{n-k}(x))
-\Delta(B_k(x+y))E_{n-k}(x)
\\=&k(x+y)^{k-1}2x^{n-k}+B_k(x+y)2x^{n-k}-k(x+y)^{k-1}E_{n-k}(x).
\endalign$$
Thus
$$\align &\Delta^*\(\sum_{k=1}^n\f{B_k(x+y)}kE_{n-k}(x)\)
\\=&2\sum_{k=1}^n(x+y)^{k-1}x^{n-k}+2\sum_{k=1}^n\f{B_k(x+y)}kx^{n-k}
-\sum_{k=1}^n(x+y)^{k-1}E_{n-k}(x)
\\=&2\f{(x+y)^n-x^n}{(x+y)-x}+2\sum_{l=1}^n\bi nl\f{B_l(y)}lx^{n-l}+2H_nx^n
\\&-\sum_{l=0}^{n-1}\bi n{l+1}E_l(-y)(x+y)^{n-1-l}
\endalign$$
where we apply Lemma 2.2 in the last step.
It follows that $\Delta^*(L(x))=\Delta^*(R(x))$ where $L(x)$ and $R(x)$ are the
left hand side and the right hand side of (2.5) respectively.
Applying Lemma 2.1(ii) we find that $L(x)=R(x)$.

The proof of Theorem 1.2 is now complete. \qed

\heading{3. Final remarks}\endheading

Quite recently Dunne and Schubert [DS] proposed the following generalization of Miki's and Matiyasevich's
identities involving
the well-known gamma function.
$$\align
&\f1{\Gamma(2n+2p)}\sum_{k=1}^{n-1}\frac{B_{2k}B_{2n-2k}}{8k(n-k)}
\cdot\frac{\Gamma(2k+p)\Gamma(2n-2k+p)}{\Gamma(2k)\Gamma(2n-2k)}\\
=&\Gamma(p+1)\sum_{k=1}^{n}\frac{B_{2k}B_{2n-2k}\Gamma(2k+p)}{(2k)!(2n-2k)!\Gamma(2k+2p+1)}
+\frac{B_{2n}}{(2n)!}\sum_{l=1}^{2n-1}\beta(l+p,p+1),
\endalign$$
where $n\in\{2,3,\ldots\}$, $p\gs0$ and
$$\beta(a,b)=\int_{0}^1 t^{a-1}(1-t)^{b-1}dt=\f{\Gamma(a)\Gamma(b)}{\Gamma(a+b)}.$$

 This result follows from our following theorem in the special case $p=q$ and $x=0$.

\proclaim{Theorem 3.1} Let $n>1$ be an integer, and let $p\gs0$ and $q\gs0$. Then
$$\aligned
&\frac{\Gamma(n)}{\Gamma(n+p+q)}
\sum_{k=1}^{n-1}B_k(x)B_{n-k}(x)
\frac{\Gamma(k+p)\Gamma(n-k+q)}{k!(n-k)!}\\
=&\sum_{l=2}^{n}\bi{n-1}{l-1}\f{B_l}lB_{n-l}(x)
\frac{\Gamma(l+p)\Gamma(q+1)+\Gamma(l+q)\Gamma(p+1)}{\Gamma(l+p+q+1)}\\
&+\f{B_{n}(x)}n(H_n(p,q)+H_n(q,p)),
\endaligned\tag 3.1$$
where
$$H_n(p,q)=\sum_{k=1}^{n-1}\beta(k+q,p+1)
=\cases \beta(p,q+1)-\beta(p,n+q)&\t{if}\ p\not=0,\\
\sum_{k=1}^{n-1}(k+q)^{-1}\quad &\t{if}\ p=0.\endcases$$
\endproclaim

 In the case $p=q=0$ or $p=q=1$ Theorem 3.1 yields (1.6) or (1.7).

We omit our proof of Theorem 3.1 which is very similar to the proof of Theorem 1.1.
The key point is the following extension of Chu's identity which can be easily proved by induction on $n$.
\proclaim{Lemma 3.1} Let $n\gs l>0$ be integers, and let
$p\gs0$ and $q>0$. Then
$$\sum_{k=l}^{n}\bi{n-l}{k-l}\beta(k+p,n-k+q)=\beta(l+p,q).\tag3.2$$
\endproclaim

\bigskip

 Department of Mathematics, Nanjing University, Nanjing 210093,
The People's Republic of China.

{\it E-mail}: {\tt haopan79\@yahoo.com.cn}

\medskip
 Department of Mathematics (and Institute of Mathematical Science),
Nanjing University, Nanjing 210093, The People's
Republic of China.

{\it E-mail}: {\tt zwsun\@nju.edu.cn}

Homepage: {\tt http://pweb.nju.edu.cn/zwsun}

\enddocument